\documentclass[12pt,twoside]{article}
\usepackage[dvips]{graphicx}
\usepackage{latexsym}
 
\usepackage{amssymb}

\def\R{{\mathbb R}}

\newtheorem{theorem}{Theorem}

\newtheorem{definition}[theorem]{Definition}

\def\R{{\mathbb R}} 
\def\C{{\mathbb C}} 

\begin{document}
\overfullrule=0pt
\baselineskip=24pt
\font\tfont= cmbx10 scaled \magstep3
\font\sfont= cmbx10 scaled \magstep2
\font\afont= cmcsc10 scaled \magstep2
\title{\tfont A Review of Alexander Isaev's\\
{\it Spherical Tube Hypersurfaces} }

\bigskip
\author{Thomas Garrity\
\\ Department of Mathematics\\ Williams College\\ Williamstown, MA  01267\\
tgarrity@williams.edu}

\date{}
\maketitle

As this book is about a subfield of CR geometry,  we should first attempt to answer,  ``What is CR geometry and why should we care?''  For CR geometry, the archetype question is: given sets $\Gamma_1$ and $\Gamma_2$ in $\C^n$, when is there a biholomorphic map $\phi:\C^n\rightarrow \C^n$ that take $\Gamma_1$ to $\Gamma_2$?  From this question stems a fascinating interplay of differential geometry, several complex variables and partial differential equations, with a lot of linear algebra thrown into the mix.

\section{Sets in $\C$} Even for sets in $\C,$ the question of when sets $\Gamma_1$ and $\Gamma_2$ are biholomorphically equivalent leads to great mathematics.

Let us start with assuming the sets  $\Gamma_1$ and $\Gamma_2$ are finite (in which case at the least they must have the same number of elements).  Here a complete solution exists.  What is key, given three distinct complex points $z_1, z_2$ and $z_3,$ is that there is a unique biholomorphic map such that:
$$z_1\rightarrow 0, z_2\rightarrow 1, z_3\rightarrow 2.$$
In fact, using M\"{o}bius transformations,  this map is explicitly
$$\phi(z) = \frac{  z(2z_3-2z_2)+ 2(z_1z_2-z_1z_3)}{z(z_1-2z_2+z_3) + (z_1z_2-2z_1z_3   +z_2z_3)}.  $$
Then if $\Gamma_1=\{z_1,z_2,z_3\}$ and $\Gamma_2=\{w_1,w_2,w_3\}$, let
$\phi_i$ send $\Gamma_i$ to $\{0,1,2\}$.  Our desired map is simply
$\phi_2^{-1}\circ \phi_1.$
(This map is not actually a biholomorphic map from $\C$ to itself but instead a map from the complex projective line to itself.)

At least as interesting is  the Riemann Mapping Theorem, which says that any two simply connected open sets, neither of which are all of $\C$. are biholomorphically equivalent. As we will see, the analogous result for $\C^n$, $n\geq 2$, is profoundly false.

Thus even for $n=1$, we are led to some significant and deep mathematics.  But the study of biholomorphic maps for $n=1$ is not usually labelled as a part of CR geometry.

\section{No Riemann Mapping  Theorem for $n\geq 2$ }
Poincaire around 1900 showed that there is no biholomorphic map taking the unit ball $B$ in $\C^2$:
$$B=\{(z_1,z_2)\in \C^2: |z_1|^2+|z_2|^2 <1\}$$
to the polydisc
$$P=\{(z_1,z_2)\in \C^2: |z_1|<1,\;|z_2| <1\}.$$
Hence there is no possibility for a Riemann Mapping Theorem in $\C^2$.  Even further, in work that leads to much of current  CR geometry, Poincaire showed that there is no biholomorphic map taking the boundary of $B$, the unit sphere $S^3$, to the  boundary of the polydisc $\partial (P).$  Note that both of these are real three dimensional submanifolds of the complex space $\C^2.$ This will lead us to the study of real hypersurfaces in $\C^n$, which, as we will see, are the model examples of CR structures.

\section{Linear Algebra Interlude}
We want to understand real subspaces of the complex vector space $\C^n$.  Let $z_1, \ldots, z_n$ be complex coordinates for $\C^n$, with each $z_k=x_k+iy_k.$  For $\C^n$, we have the natural linear map of multiplication by $i$:
$$i (z_1,\ldots, z_n)=(iz_1,\ldots, iz_n).$$

We can identify the complex $n$-dimensional space $\C^n$ with the real $2n$-dimensional vector space
$ \R^{2n}$ by  the map:
$$  (z_1,\ldots, z_n)\rightarrow  (x_1,y_1,\ldots, x_n,y_n).$$
The map $i$ will correspond to a linear map $J:\R^{2n}\rightarrow \R^{2n}$ so that the following diagram is commutative:
$$\begin{array}{ccc}
\C^n & \overrightarrow{i} & \C^n \\
\downarrow && \downarrow \\
\R^{2n} & \overrightarrow{J} & \R^{2n}.
\end{array}$$
We have
$$J(x_1,y_1,\ldots, x_n,y_n) = (-y_1, x_1, \ldots , -y_n, x_n).$$
Thus,  $J$  is the block diagonal matrix with diagonal entries the two-by-two matrices
$$\left(\begin{array}{rr} 0 & -1 \\ 1&0 \end{array}\right),$$
giving us that  $J^2=-I$.
We are interested in the interplay between real subspaces of $\C^n$  and the complex structure of $\C^n$.  For example, let $V$ be a real three-dimensional subspace of $\C^2.$ Thinking of $\C^2$ as $\R^4$, we know that
$J(V)$ is another real three-dimensional subspace.   Then purely for dimensional reasons we know that
$$V\cap J(V)$$
is real two-dimensional but, more importantly, that
$$J(V\cap J(V)) = V\cap J(V),$$
meaning that  $V\cap J(V)$ can be identified with a complex one-dimensional subspace of $\C^2$.

In fact, the key for us is that if $V$ is a real $(2n-1)$-dimensional subspace of $\C^n$, then  the subspace $$H=V\cap J(V)$$
will not only  be of real dimension $(2n-2)$ but more importantly be a complex subspace of complex dimension $n-1.$
Further this means that every real $(2n-1)$-dimensional subspace $V$ has a $(2n-2)$-dimensional subspace $H$
that carries a complex structure and a real subspace $W$ of real dimension one so that
$$V=H\oplus W.$$
The subspace $H$ is canonical.  The real subspace $W$ is not.

\section{ CR structures}

Suppose $M$ is a smooth real hypersurface in $\C^n$.  For example, $M$ might be the unit sphere
$$S^3=\{(z_1,z_2)\in \C^2:|z|^2 + |z_2|^2=1\}.$$
Of course, identifying $\C^n$ with $\R^{2n}$, we can think of $M$ as a smooth real $(2n-1)$-dimensional submanifold of $\R^{2n},$ which for the unit sphere would mean writing $S^3$ as
$$S^3=\{(x_1,y_1,x_2,y_2)\in \R^4: x_1^2+ y_1^2 + x_2^2 + y_2^2 =1\}.$$

By smoothness, we know that at every point $p\in M$, there is the real $(2n-1)$-dimensional tangent space
$T_pM.$  But we do not want to forget that $M$ is lying in a complex space.  From our linear algebra interlude, we have that  at every point $p\in M$, there is  a real $(2n-2)$-dimensional subspace $H_p$  of $T_pM$:
$$H_p=T_pM\cap J(T_pM).$$
As we have seen, $H_p$ can be thought of as the part of the vector space $T_pM$ that inherits a complex structure.  Hence at each point $p\in M$  we can write the tangent space as a direct sum:
$$T_pM=H_p\oplus W_p,$$
where $H_p$ is a a real $2n-2$ vector space that inherits a complex structure and hence can be identified with a copy of $\C^{n-1}$ and $W_p$ is a real one-dimensional vector space and hence  has no complex structure.

All of this can be extended to a splitting of the tangent bundle $TM=H\oplus W$.  Here we have
$$J:H\rightarrow H.$$
Since $J^2=-I$, we have $\C\otimes H$ splitting into two subbundles: $\mathcal{H}$ the complex rank $n-1$ subbundle corresponding to the $i$-eigenspace of $J$ and its complex conjugate  $\overline{\mathcal{H}}$, the subbundle corresponding to the $-i$ eigenspace of $J$.  If we write $\mathcal{W}=\C\otimes W$, we have
$$\C\otimes TM = \mathcal{H}\oplus \overline{\mathcal{H}}\oplus \mathcal{W}.$$
All of this can be done to any smooth real hypersurface in $\C^n$.

A CR structure is the generalization of the above to abstract manifolds:
\begin{definition} Let $M$ be a real differentiable manifold of dimension $d$.  Suppose that $2n+c=d$.  The manifold $M$ is a {\rm CR structure of codimension $c$} if there is  a complex rank $n$ subbundle
$$\mathcal{H} $$
of the complexified tangent bundle $\C \otimes TM$ such that
\begin{enumerate}
\item $\mathcal{H} \cap \overline{\mathcal{H}}$ is the zero subbundle
\item $[\mathcal{H}, \mathcal{H}]  \subset \mathcal{H}$.
\end{enumerate}
\end{definition}
(Here $[\mathcal{H}, \mathcal{H}]  \subset \mathcal{H}$ means that  for all local smooth sections $v$ and $w$ of $\mathcal{H}$ considered as vector fields, the commutator $[v,w]=v\circ w - w\circ v$ is also a smooth section of $\mathcal{H}.$)

In particular, hypersurfaces are CR structures of codimension one.   (Also, the condition $[\mathcal{H}, \mathcal{H}] \subset \mathcal{H}$ is always satisfied  for a real hypersurface in a complex vector space; this condition becomes more relevant when we ask various types of embedding questions in a few sections.)
 We now want to find a natural notion of equivalence between two CR structures.  This equivalence should reduce to the idea of biholomorphism of the ambient complex space, which means that not only should the two structures be  diffeomorphic but also that the complex parts of the tangent bundles should map to each other.  Let us make this into a rigorous definition.  Let $M_1$ be a CR structure of codimension $c$ with codimension $c$ subbundle $H_1$.  Further, let the involution corresponding to the complex structure be denoted as
 $$J_1:H_1\rightarrow H_1.$$
Let $M_2$ denote another CR structure of codimension $c$ with corresponding subbundle $H_2$ with involution $J_2$.
 \begin{definition} The CR structures $M_1$ and $M_2$ are {\rm CR equivalent} if there is a diffeomorphism
 $\pi:M_1\rightarrow M_2$
 such that
 \begin{enumerate}
 \item $\pi_*(H_1)\subset H_2$
 \item $J_2\circ \pi_* = \pi_* \circ J_1.$
 \end{enumerate}
 \end{definition}

 Much of the study of CR geometry is an attempt to determine when two CR structures are equivalent.

 \section{ The Levi Form}
 The most important tool for studying CR structures is the Levi form, which plays the role of curvature in CR geometry.  (More accurately, the Levi form is the CR analog of the second fundamental form of classical differential geometry.)  As with different types of curvature  in differential geometry, there are a number of ways for defining Levi forms, each with its  strengths and weaknesses.  Also, as is the case for curvature, the machinery can quickly become difficult and abstract.

 We start with the case  when our CR structure $M$ is a real smooth hypersurface of $\C^{n+1}$.  The Levi form will attach to each point $p\in M$  an $n\times n$ Hermitian matrix $L=(L_{ij})$.  We will first show how to compute this matrix and only then give some justification for its meaning.

  Since $M$ is an embedded smooth hypersurface,  for any point $p\in M$, there is a smooth function $\rho$ so that near $p$,  $M$ is given by
 $$\{(z_1, \ldots , z_{n+1})\in \C^{n+1}:\rho(z_1, \ldots , z_{n+1})=0 \;\mbox{and}\; \mathrm{d}\rho(z_1, \ldots , z_{n+1})\neq 0\}.$$
 Form the $(n+1)\times (n+1)$ Hermitian matrix of second partial derivatives:
  $$\tilde{L} = \left( \frac{\partial^2 \rho}{\partial z_i \partial \overline{z}_j}\right).$$
 The Levi form  $L$ is the restriction of $ \tilde{L} $ to the subbundle $\mathcal{H}$.  More precisely, for each
  $p\in M$, choose a basis $w_1, w_2, \ldots , w_n$ for the complex vector space $\mathcal{H}_p.$ We can write each of these basis elements as
  $$w_i=\sum_{k=1}^{n+1} a_{ik} \frac{\partial }{\partial z_k}.$$
  Then
  basis vectors for the vector space $\overline{\mathcal{H}}_p$ are the conjugates $\overline{w}_1, \ldots , \overline{w}_n$.  The entries of the Levi form are:
 $$L_{ij} = \sum_{k,l=1}^n \frac{\partial^2 \rho}{\partial z_k \partial \overline{z}_l} a_{ik}\overline{a}_{jl}.$$
 This is just to show that the Levi form can be easily calculated.

  All of this depends on the  choices  for the  defining function $\rho$,  the coordinates $z_1, \ldots , z_{n+1}$ for  $\C^{n+1}$  and the  basis $w_1, w_2, \ldots , w_n$ for  $\mathcal{H}.$   Change any of these  and the Levi form changes.  Luckily though, for  any of these possible different choices, there will be an invertible Hermitian $n\times n$ matrix $A$ such that the matrix $L$ becomes $\overline{A}^T L A$.  It can be shown that if two CR structures are equivalent, then the Levi forms must  map to each other.   Thus the standard invariant theory for Hermitian matrices can be used to understand the CR equivalence problem.  For example, the most studied real hypersurfaces are those for which the Levi form is positive definite; such hypersurfaces are called {\it strongly pseudoconvex}.

  Now let us see why the Levi form is, to some extent, a natural object to study.  First, it certainly is not uncommon to have a manifold given as the zero locus of a function $\rho$ and form the matrix of its second derivatives.  If $t_1, \ldots , t_{2n+2}$ are real coordinates for $\C^{n+1}$ treated as a real vector space,  then the second fundamental form (also called the Hessian) is the $(2n-1)\times (2n-1)$ symmetric matrix obtained by restricting the $2n\times 2n$ matrix
  $$\left( \frac{\partial^2 \rho}{\partial t_i \partial t_j}   \right)$$
  to the tangent space $TM.$  The Levi form can be thought of as the complex analog of the second fundamental form.  But this link is far stronger than mere analog.  As shown, for example, in  \cite{Krantz01}:
  $$\mbox{Second Fundamental Form} = \mbox{Levi Form} +\mbox{other stuff}.$$
  The second fundamental form contains all curvature information about $M$ (see for example Morgan's {\it Riemannian Geometry: A Beginners Guide} \cite{Morgan98}).

  We want to know when  two real hypersurfaces of $\C^n$ can be mapped to each other by a biholomorphic map.  Standard differential geometry is of little help, as curvature properties are not preserved under biholomorphic maps.  The Levi form, though,  can be shown to be precisely the part of the second fundamental form that transforms decently under biholomorphic maps.  It is the closest we can get to coupling curvature ideas with CR equivalences.

  But all of the above is for an embedded hypersurface.   With enough work, one can show that the following definition, which makes sense for any CR structure, will agree with above calculations:
   \begin{definition} The {\rm Levi form} is the bundle map
   $$L:\mathcal{H} \times \mathcal{H} \rightarrow \C \otimes TM/\left( \mathcal{H}\oplus \overline{\mathcal{H}} \right)$$
   defined by
   $$L(X,Y) = i\pi \left[\overline{X}, Y\right],$$
   where $X$ and $Y$ are sections of $\mathcal{H}$  and where $\pi:\C\otimes TM \rightarrow  \C \otimes TM/\left( \mathcal{H}\oplus \overline{\mathcal{H}} \right)$ is the natural projection.
\end{definition}
The best way to see that this agrees with our first definition is to work out the details for  the unit sphere $S^3$ in $\C^2$.

Thus the Levi form is a vector-valued Hermitian form.  Classical invariant theory can now be applied to trying to understand the equivalence problem of CR structures \cite{Garrity-Mizner} \cite{Garrity-Grossman} \cite{Mizner}  \cite{Garrity-Mizner97}.

 \section{Three motivations for CR geometry}
 We have been motivating the study of CR geometry via asking for when two CR structures are CR equivalent, which is a natural generalization of trying to answer the  geometric question of when two real submanifolds of $\C^n$ are equivalent under a biholomorphic map.  There are two other areas of mathematics that also naturally lead to  CR geometry.

 The first is the study of domains of holomorphy.   Complex analysis in one dimension is a fundamentally different subject than complex analysis in several dimensions.  For example,  as mentioned earlier, the Riemann Mapping Theorem is  true only for  dimension one.  But there are other profound differences.   Consider the one-variable function $f(z)=1/z.$ This function $f(z)$ is holomorphic on the punctured unit disc $\{z\in \C:0<|z|^2 <1\}.$  This function cannot be extended to a holomorphic function on the entire unit disc.  There is a true pole at the origin.  A natural question  to ask is if there is a   function $f(z_1,z_2)$ that is holomorphic on the punctured ball
 $\{(z_1,z_2)\in \C^2: 0<|z_1^2 +|z_2|^2 < 1\} $ that cannot be extended to the entire ball.  In other words, is there a function that is holomorphic in two-variables that has a pole at an isolated point?  Certainly there are smooth  functions with isolated poles  (such as $1/(|z_1|^2+|z_2|^2))$.  The punchline of Hartogs theorem:
 \begin{theorem}  Let $U$ be an open connected region in $\C^n$, with $n>1$ and let $V$ be a compact connected region contained in $U$.  Then any function $f(z_1,\ldots, z_n)$ that is holomorphic on $U-V$ can be extended to a function that is holomorphic on all of $U$.
 \end{theorem}
\noindent is that such isolated poles do not exist.  This leads to:
\begin{definition} A domain $D$ in $\C^n$ is a {\rm domain of holomorphy} if for any point $p\in \partial (D)$, there is a function $f(z_1,\ldots ,z_n)$ that is holomorphic on $D$ that cannot be extended to $p$.
\end{definition}
 Thus the punctured ball in $\C^2$ is not a domain of holomorphy.  The ball is, though.  Let $p=(a,b)\in \partial (D)$.  The function
 $$f(z_1,z_2) = \frac{1}{\overline{a} z_1 + \overline{b}z_2- 1}$$
 will be holomorphic in the ball but has a pole at $p$.

 The question becomes trying to find a reasonable condition to identify domains of holomorphy.  Subject to certain smoothness restraints, the boundary of any domain $D$ is a smooth real hypersurface in $\C^n.$  The quite spectacular fact is that such $D$ will be a domain of holomorphy when the Levi form of the boundary is positive definite (in other words, when every boundary point is strongly pseudoconvex).  Hence the CR geometry of the hypersurface tells us about the domain.

 The second  is the study of existence (more precisely, the failure of existence) of solutions of  some linear partial differential equations.  By the early 1950s, it was shown that there is always a solution $f(t_1, \ldots, t_n)$
 to the linear partial differential equation
 $$a_1\frac{\partial f}{\partial t_1} +\cdots + a_n\frac{\partial f}{\partial t_n} =g(t_1, \ldots, t_n)$$
 where the  $a_i$ are constants and $g$ is any smooth function.  At the time, it was believed that any reasonable linear partial differential equation would always have solutions.  It thus came as quite a shock in 1955 when Hans Lewy showed
 that  a solution 
 $f(t_1, t_2,t_3)$ to
 $$\frac{\partial f}{\partial t_1} +i\frac{\partial f}{\partial t_2}  -2i(x+iy)\frac{\partial f}{\partial t_3} =g(t_1, t_2,t_3)$$
 will exist 
 only under the drastically restrictive assumption that $g$ is real-analytic.  Though the proof is not hard, there was no real indication as to why this particular PDE had no solutions for $g$ that were merely smooth.  Luckily, in 1974, Nirenberg showed that the failure of the above Lewy PDE to have a solution corresponded to the existence of a smooth codimension one, three-dimensional CR structure that was not CR equivalent to a hypersurface in $\C^2$.  Suddenly, people saw geometric obstructions for solving the Lewy PDE.  This in turn has led to a rich interplay between the study of solutions to linear PDEs and embeddings of CR structures.  For example, Kuranishi showed that all codimension one strongly pseudoconvex CR structures of dimension nine or greater can be realized as hypersurfaces in complex space.  Akahori non-trivially extended this to the dimension seven case.  One of the big open questions is determining what happens for dimension five.

 \section{Getting Started in CR Geometry}

 While the first three chapters  of  {\it Spherical Tube Hypersurfaces} do  cover the basics of CR geometry, it is probably not the place to start for a novice.  Luckily, a number of more introductory texts have been written in the last twenty or so years.  Reflecting the various sources that have shaped CR geometry, each of these books has a different feel and emphasis.  Jacobowitz's {\it An Introduction to CR Structures}  \cite{Jacobowitz90}  is probably the best text to get you as quickly as possible to the heart of this book, since he does a great job developing the Chern-Moser machinery, which is critical for  much of Isaev's work.

  Boggess's {\it CR Manifolds and the Tangential Cauchy-Riemann Complex} \cite{Boggess91} takes quite a different approach.  In Part I, Boggess does an excellent job of developing the basics of multi-variable function theory.  In Part II, he presents the basics of CR structures.  The last two parts concern the function theory of CR structures, in the following sense.  With coordinates $z_1, \ldots z_n$ for $\C^n$, consider the operator
  $$\overline{\partial}f = \sum_{k=1}^n \frac{\partial f}{\overline{\partial}z_k} \mathrm{d}\overline{z}_k,$$
  where, if $z_k=x_k+iy_k,$
  $$\frac{\partial f}{\overline{\partial}z_k} = \frac{1}{2}\left(\frac{\partial f}{\partial x_k} + i\frac{\partial f}{\partial y_k}  \right).$$
  The classical fact is that a smooth function $f$ on $\C^n$ is complex analytic if and only if $\overline{\partial}f=0.$
  For a CR structure $M$ there is an analogous operator $\overline{\partial}_M$.  The dream (which is not true) is that,  for any smooth function $f$ defined on a real CR manifold  $M$ in $\C^n,$ if  $\overline{\partial}_M(f)=0$, then $f$ is the restriction of a function that is holomorphic in a neighborhood of $M$.  This dream is true, though, if $M$ is real-analytic and if the function $f$ is real-analytic.  The operator $\overline{\partial}_M$ can be defined for any CR structure.  A function $f$ on $M$ is said to be a  ${\it CR}$ ${\it  function}$ if  $\overline{\partial}_M(f)=0$.  CR functions  are suppose to mimic those functions that are the restrictions of holomorphic functions.  Part III of Boggess talks about when CR functions for $M$ in $\C^n$ are actual restrictions of holomorphic functions.  Part  IV concerns the tangential Cauchy-Riemann complex of the title.  This complex arises in a natural fashion from the operator $\overline{\partial}_M$.

 Pitched to an audience at a  slightly higher higher level of mathematical maturity is the wonderful {\it Real Submanifolds in Complex Space and Their Mappings} by Baouendi, Ebenfelt and Rothschild \cite{Baouendi-Ebenfelt-Rothschild99}.  In fact, the preface accurately  states that the ``material in this book is intended  to be accessible to mature graduate students; no previous knowledge of several complex variables is assumed of the reader.''  As the tile suggests, this book concentrates on embedded CR structures and on the great deal of work done in the 1980s and 1990s,  rich work that is ongoing.

 As its title indicates,  {\it Differential Geometry and Analysis on CR Manifolds} by Dragomir and Tomassini \cite{Dragomir-Tomassini06} emphasizes the differential geometry of CR structures.  Traditional differential geometry is a source of many interesting systems of partial differential equations.  The original geometric inspiration for each of these systems can then be used to understand the solutions of these systems.   The CR curvature analogs give rise to similar PDE systems.  As just one example, Chapter 3 of this book deals with the CR-Yamabe problem.   In differential geometry, the Yamabe equation corresponds to trying to solve the following geometric problem.  Let $M$ be a Riemannian manifold with a fixed metric $g$.   Is there another metric $g'$ on $M$ that is conformal with respect to $g$ and, most importantly, has constant scalar curvature?   In \cite{Webster78}, Webster constructed a CR analog to the Riemann curvature tensor and hence a CR analog of the Riemannian scalar curvature. Thus there  is a CR analog of the Yamabe problem, leading to a CR-Yamabe equation.  These types of analogs are not obvious, and more importantly, attempts to find solutions are rarely easy and have led (and will continue to lead) to a lot of great mathematics.

 D'Angelo's {\it Several Complex Variables and the Geometry of Real Hypersurfaces} \cite{D'Angelo93} concentrates on what are called higher order invariants of real hypersurfaces.   Though definitely a book on analysis, it has more of an algebraic feel (in particular a commutative algebraic feel)  than the other books.  In fact, D'Angelo states in the Preface, `` I believe that combined use of the algebraic and analytic ideas ... forms a useful tool for attacking several of the field's open problems.  In particular consider the following analogy.  Strongly pseudoconvex points correspond to the maximal ideal, while points of finite type correspond to ideals primary to the maximal ideal.  Making sense of this simple heuristic idea is perhaps the raison d'\^{e}tre for the writing of this book.''

 Huang's survey paper  `Local Equivalence Problems for Real Submanifolds in Complex Space' \cite{Huang04} provides a good overview of CR equivalence problems.  There are also the three great recent  surveys in the {\it Notices of the American Mathematical Society}:  Treves's  `A Treasure Trove of Geometry and Analysis:The Hyperquadric' \cite{Treves00}, D'Angelo and Tyson's `Cauchy-Riemann and Sub-Riemannian Geometries' \cite{D'Angelo-Tyson10} and  Ezhov, McLAughlin and Schmalz's `From  Cartan to Tanaka: Getting Real in the Complex World' \cite{Ezhov-McLAughlin-Schmalz11} .

 \section{Spherical Hypersurfaces are CR Flat}
 The book under review is overwhelmingly  concerned with the equivalence problem for CR structures.  What has happened over the years is that people have discovered different natural classes of CR structures, and then concentrate on the structure of each of these classes.  Thus this text is concerned with those CR structures that are hypersurfaces in some $\C^n$ that are simultaneously spherical and tube hypersurfaces, both of which we will define in a moment.

 For now though, we approach the equivalence problem in a more general fashion.  E. Cartan started this process for three dimensional CR structures in $\C^2$ by applying his method of moving frames (for an introduction to moving frames, see \cite{Ivey-Landsberg03}).  This led Chern and Moser in the 1960s to deeply extend this work, developing what is now called Chern-Moser theory.  An excellent overview of this work is in the already mentioned  {\it Notices} article  by Ezhov, McLAughlin and Schmalz \cite{Ezhov-McLAughlin-Schmalz11}.

 In traditional differential geometry, a manifold $M$ in $\R^n$ should be flat if $M$ is linear.  The various notions of curvature for more general $M$ in $\R^n$ are attempts to measure how far  $M$ is from a linear space.   In the world of CR structures, coming from Chern-Moser, the CR analog of linear (or flat) for a hypersurface $M$ in $\C^{n+1}$ is when $M$ can locally  be put into the form
 $$\{(z_1, \ldots, z_{n+1})\in \C^{n+1}: z_{n+1}+\overline{z}_{n+1}= \sum_{i,j=1}^n a_{ij}z_i \overline{z_j}\}$$
 where the matrix  $A= (a_{ij})$ of constants is Hermitian.  The quintessential example is an actual sphere.
 The Chern-Moser curvarture is a measurement of how far a CR structure is from one that is flat.  Thus CR-flat should mean that $M$ is a quadric hypersurface of a Hermitian matrix.  These are also called spherical hypersurfaces, leading to the somewhat perverse terminology that a CR hypersurface is flat when it is spherical.  The term ``flatness'' comes from the curvature interpretation while the term ``spherical" comes from the quadric surface interpretation.

 \section{Tube Domains}
 The other condition  that Isaev places on CR hypersurfaces is to require them to be tube domains. The study of tube domains for spherical CR structures was pioneered by  Yang \cite{Yang82}.   A CR hypersurface $M$ in  $\C^{n+1}$ will be a tube domain if the following happens.  There must be a totally real subspace $V$  of real dimension $n+1$ in $\C^{n+1}$, which means that
 $$V\cap iV=(0)$$
  and a real $n$ dimensional hypersurface $M_{\R}$ in $V$ such that
  $$M=M_{\R} \times iV.$$
  $M_{\R}$ is called the {\it base}.

  For example, in $\C^3$,  the hypersurface
  $$M=\{(z_1, z_2,z_3)\in \C^3:  \mbox{Re}(z_3)= (\mbox{Re}(z_1))^2 + (\mbox{Re}(z_2))^2  \}$$
  is a tube domain; simply let $V=\mbox{Span}( x_1, x_2, x_3)$ and set
  $$M_{\R}= \{ (x_1,x_2,x_3)\in V: x_3 = x_1^2+x_2^2  \}.$$

 \section{The Invariant Theory of Spherical Tube Hypersurfaces}
 The main goal and purpose of Isaev's  book  is to explore the invariant theory of the special class of spherical tube hypersurfaces.  To be clear, this is a major restriction on  the generality of CR hypersurfaces, but a restriction that results in a manageable classification theory.  This type of move is of course standard in mathematics.  After all, for example, no one would expect an easy classification theorem of all topological spaces.  Instead, people make restrictions, such as the classical work of looking only at compact surfaces.

 After spending the first three chapters on the necessary foundations of CR structures, in chapter four Isaev shows that   spherical tube hypersurfaces reduce to three broad classes.   This classification has  a heavy linear algebraic feel, which is a direct consequence of the spherical and tube  conditions.

 Here is a hint of how to proceed.  If $M$ is a spherical hypersurface in  $\C^{n+1}$, then $M$ can be written as  $\{(z_0, \ldots, z_{n})\in \C^{n+1}: z_{0}+\overline{z}_{0}= \sum_{i,j=1}^n a_{ij}z_i \overline{z_j}\}$.  Then $M$ is $(k, n-k)$ spherical if the Hermitian matrix $(a_{ij})$ has signature $(k,n-k)$.  If $M$ is also a tube domain, then  we can write $M$ as $M=M_{\R} \times iV.$  From section 3.1 of the book, 
 
 \begin{quote}There is a natural equivalence relation for tube hypersurfaces. Namely, two tube hypersurfaces $M_1, M_1$ are called {\it affinely equivalent} if their bases are affinely equivalent in $\R^{n+1}$, i.e. if there exists an affine transformation of $\C^{n+1}$ of the form 
 $$Z\mapsto CZ + b, \;\; C\in GL(n+1, \R), b\in \R^{n+1}$$
 that maps $M_1$ to $M_2$.... If $M$ is a locally closed tube hypersurface and $p\in M_{\R}$,  then there exist a tube hypersurface $\mathcal{M}$ of $p$  and an affine transformation  $\mathcal{A}$ of $\C^{n+1}$  ... such that (i)  $ \mathcal{A}(p) = 0$,  (ii) $\mathcal{A}(\mathcal{M}) = \Pi^{-1}(\mathcal{V})$ for a neighborhood $\mathcal{V}$ of the origin in $\R^{n+1}$ of the form $\mathcal{V}= \Omega + I$, where $\Omega$ is a domain in $\R^n$ and $I$ is an interval in the line $\{x=0\}$, with $\R^n$ identified with the linear subspace of $\R^{n+1}$ given by $x_0=0$, and (iii) for $W:= M\cap \mathcal{M}$ the base of the tube hypersurface $\mathcal{A}(W)$ is represented in $\mathcal{V}$ as a graph 
 $$x_0=F(x) \;\mbox{with} \; F(0)=0, \frac{\partial F}{\partial x_{\alpha}} (0) =0,$$
 where $F$ is a function on $\Omega.$ \end{quote}
Thus we are trying to put the original spherical tube hypersurface $M$ into somewhat of a more canonical form.   This form for $M$ is called a {\it standard representation}.  Note that we  cannot call this {\it the} standard representation, as it is not unique.  

One would like to see what restrictions and structures can be placed on the now important function $F$.   For a spherical tube hypersurface $M$ of signature $(k,n-k)$, we can put  $M$ into a  standard representation such that 
the function $F$ satisfies 

$$\frac{\partial^2 F}{\partial x_{\alpha} \partial x_{\beta}} = \sum_{\gamma = 1}^{n+1} \frac{\partial F}{\partial x_{\gamma}}\left( D_{\alpha}^{\gamma} \frac{\partial F}{\partial x_{\beta}} +  D_{\beta}^{\gamma} \frac{\partial F}{\partial x_{\alpha}}  +C_{\alpha \beta}^{\gamma}        \right) + H_{\alpha \beta}, $$
where $D_{\alpha}^{\gamma}, C_{\alpha \beta}^{\gamma} $ and $H_{\alpha \beta}$ are real constants.  (Here the $x_{\alpha}$ are the real parts of the complex coordinates $z_{\alpha} = x_{\alpha} + iy_{\alpha}$ of the ambient $\C^{n+1},$ for $\alpha = 1, \ldots ,n$.)  Further, the $H_{\alpha \beta}$ can be chosen so that 

$$H_{\alpha \beta} = \left\{   \begin{array}{cl} 1 & \alpha=\beta, \alpha =1, \ldots ,k \\
									-1 &  \alpha=\beta, \alpha =k+1, \ldots ,n  \\
									0 & \alpha \neq \beta \end{array}   \right.$$

									Isaev uses these constants to construct his three classes of spherical tube hypersurfaces.  For the rest of the book he explores the consequences of placing natural restrictions on these constants, resulting in many cases to complete classifications, at times into a finite number of possibilities and at times into infinite families.
									
 As can be imagined, the techniques and results quickly become, by necessity,  quite technical. These are also quite interesting.    Much of this is the work   of the author and also of Fels and Kaup (and a number of other people).  
 
 Thus this book will be of interest and of value to everyone working  on the equivalence problem for CR structures.

\bibliographystyle{amsplain}

\begin{thebibliography}{99}

\bibitem{Baouendi-Ebenfelt-Rothschild99} M.  Baouendi, P. Ebenfelt and L. Rothschild, {\it Real Submanifolds in Complex Space and Their Mappings}, Princeton Mathematical Series, Number 47, Princeton University Press, 1999.



\bibitem{Boggess91} A. Boggess, {\it CR Manifolds and the Tangential Cauchy-Riemann Complex}, Studies in Advanced Mathematics, CRC Press, 1991.

\bibitem{D'Angelo93} J. D'Angelo, {\it Several Complex Variables and the Geometry of Real Hypersurfaces}, Studies in Advanced Mathematics, CRC Press, 1993.

\bibitem{D'Angelo-Tyson10} J. D'Angelo and J. Tyson, Cauchy-Riemann and Sub-Riemannian Geometries,  {\it Notices of the American Mathematical Society}, Vol. 57, no. 2 (2010), pp. 208-219.


\bibitem{Dragomir-Tomassini06} S. Dragomir and G. Tomassini, {\it Differential Geometry and Analysis on CR Manifolds}, Progress in Mathematics, Vol. 246, Birkh\"{a}user, 2006.

\bibitem{Ezhov-McLAughlin-Schmalz11} V. Ezhov, B. McLAughlin and G. Schmalz, From  Cartan to Tanaka: Getting Real in the Complex World, {\it Notices of the American Mathematical Society}, Vol. 58, no. 1 (2011), pp. 20-28.

\bibitem{Garrity-Grossman} T. Garrity and Z. Grossman, On Relations of Invariants for Vector-Valued Forms, {\it Electronic 
          Journal of Linear Algebra}, (2004), Vol. 11, pp. 24-40.


\bibitem{Garrity-Mizner}  T. Garrity and R. Mizner,  Invariants of Vector-Valued Bilinear and Sesquilinear Forms,  	{\it Linear Algebra and its Applications} , (1995), v. 218, pp. 225-237.

\bibitem{Garrity-Mizner97}     T. Garrity and R. Mizner,    The Equivalence Problem for Higher-Codimensional CR Structures,  	{\it  Pacific J. of Math.}, (1997), v. 177, no. 2,  pp. 211-235.	

\bibitem{Huang04} X. Huang, Local Equivalence Problems for Real Submanifolds in Complex Space, in {\it Real Methods in Complex and CR Geometry} (Lecture Notes in Mathematics, Volume 1848), Springer-Verlag, 2004, pp.109-163.


\bibitem{Ivey-Landsberg03} T. Ivey and J. M. Landsberg, {\it Cartan for Beginners: Differential Geometry Via Moving Frames and Exterior Differential Systems}, Graduate Studies in Mathematics, volume 61, American Mathematical Society, 2003.
	
	
\bibitem{Jacobowitz90} H. Jacobowitz, {\it An Introduction to CR Structures}. Mathematical Surveys and Monographs, Number 32, American Mathematical Society, 1990.

\bibitem{Krantz01}  S. Krantz, {\it Function Theory of Several Complex Variables}(AMS Chelsea Publishing), 2nd edition, American Mathematical Society, 2001.

\bibitem{Mizner}  R. Mizner, CR structures of codimension 2, {\it Journal of Differential Geometry}, Volume 30, Number 1 (1989), 167-190.

\bibitem{Morgan98} F. Morgan, \it{Riemannian Geometry: A Beginners Guide}, A K Peters/CRC Press, 1998.


\bibitem{Treves00}  F. Treves, A Treasure Trove of Geometry and Analysis:The Hyperquadric, {\it Notices of the American Mathematical Society}, Vol. 47, no. 11 (2000), pp. 1246-1256.


\bibitem{Webster78}  S. Webster,  Pseudohermitian structures on a real hypersurface, {\it Journal of Differential Geometry}, 13, (1978), pp. 25-41.

\bibitem{Yang82} P. Yang, Automorphisms of Tube Domains, {\it American Journal of Mathematics}, 104 (1982), pp. 1005-1024.
				

	





\end{thebibliography}

\end{document}